%% file: conrecreorg.tex
\documentclass[twoside,leqno]{article} \input preamble.tex

\usepackage{accents}
\newcommand{\overcirc}[1]{\accentset{\circ}{#1}}
\DeclareMathOperator{\e}{e}
\numberwithin{equation}{section}

\numberwithin{equation}{section}

\begin{document} \title{\bf Conjugate Reciprocal Polynomials with all
Roots on the Unit Circle} \author{\sc Kathleen L. Petersen \\ \sc
Christopher D. Sinclair}
\maketitle

\begin{abstract}
  We study the geometry, topology and Lebesgue measure of the set of
  monic conjugate reciprocal polynomials of fixed degree with all
  roots on the unit circle.  The set of such polynomials of degree N
  is naturally associated to a subset of $\R^{N-1}$. We calculate
  the volume of this set, prove the set is homeomorphic to the $N-1$
  ball and that its isometry group is isomorphic to the dihedral
  group of order $2N$.
\end{abstract}

\section{Introduction} 
Let $N$ be a positive integer and suppose $f(x)$ is a polynomial in
$\C[x]$ of degree $N$.  If $f$ satisfies the identity,
\begin{equation}
\label{eq:2} f(x) = x^N \overline{f(1/\overline{x})},
\end{equation} 
then $f$ is said to be {\it conjugate reciprocal}, or simply {\it CR}.
Furthermore, if $f$ is given by
\[ f(x) = x^N + \sum_{n=1}^N c_n x^{N-n}.
\] 
then (\ref{eq:2}) implies that $c_N = 1, c_{N-n} = \overline{c_n}$ for
$1 \leq n \leq N-1$ and, if $\alpha$ is a zero of $f$ then so too is
$1/\overline{\alpha}$. The purpose of this manuscript is to study the
set of CR polynomials with all roots on the unit circle.  The
interplay between the symmetry condition on the coefficients and the
symmetry of the roots allows for a number of interesting theorems
about the geometry, topology and Lebesgue measure of this set.

CR polynomials have various names in the literature including
reciprocal, self-reciprocal and self-inversive (though we reserve the
term {\it reciprocal} for polynomials which satisfy an identity akin
to (\ref{eq:2}) except without both instances of complex conjugation).

The condition on the coefficients of a conjugate reciprocal polynomial
allows us to identify the set of CR polynomials with $\R^{N-1}$.  To
be explicit, let $X_N$ be the $N -1 \times N - 1$ matrix whose $j,k$
entry is given by,
\begin{equation}
\label{eq:8}
X_N[j,k] = \left\{
  \begin{array}{ll} \frac{\sqrt{2}}{2}\left(\delta_{j,k} + 
      \delta_{N-j, k}\right) & \mbox{if} \quad 1 \leq j < \lfloor N/2 \rfloor \\ & \\
      \delta_{j,k} & \mbox{if} \quad j = N/2 
      \\ & \\ \frac{\sqrt{2}}{2}\left(i \delta_{N-j,k} - i 
        \delta_{j, k}\right) & \mbox{if} \quad      \lfloor N/2 \rfloor < j < N, \\
\end{array} \right.
\end{equation}
where $\delta_{j,k} = 1$ if $j=k$ and is zero otherwise.  For
instance, 
\[
X_5 = \frac{\sqrt{2}}{2} \begin{bmatrix}
1 & 0 & 0 & i \\
0 & 1 & i & 0 \\
0 & 1 & -i & 0 \\
1 & 0 & 0 & -i
\end{bmatrix}
\qquad \mbox{and} \qquad
X_6 = \frac{\sqrt{2}}{2} \begin{bmatrix}
1 & 0 & 0 & 0 & i \\
0 & 1 & 0 & i & 0 \\
0 & 0 & \sqrt{2} & 0 & 0\\
0 & 1 & 0 & -i & 0 \\
1 & 0 & 0 & 0 & -i \\
\end{bmatrix}.
\]
The $\sqrt{2}/2$ factor is a normalization so that $|\det X_N| = 1,$
a fact which is easily checked by induction on $N$ (odd and even cases
treated separately).

\begin{figure}[h!]
\centering
\includegraphics[scale=.6]{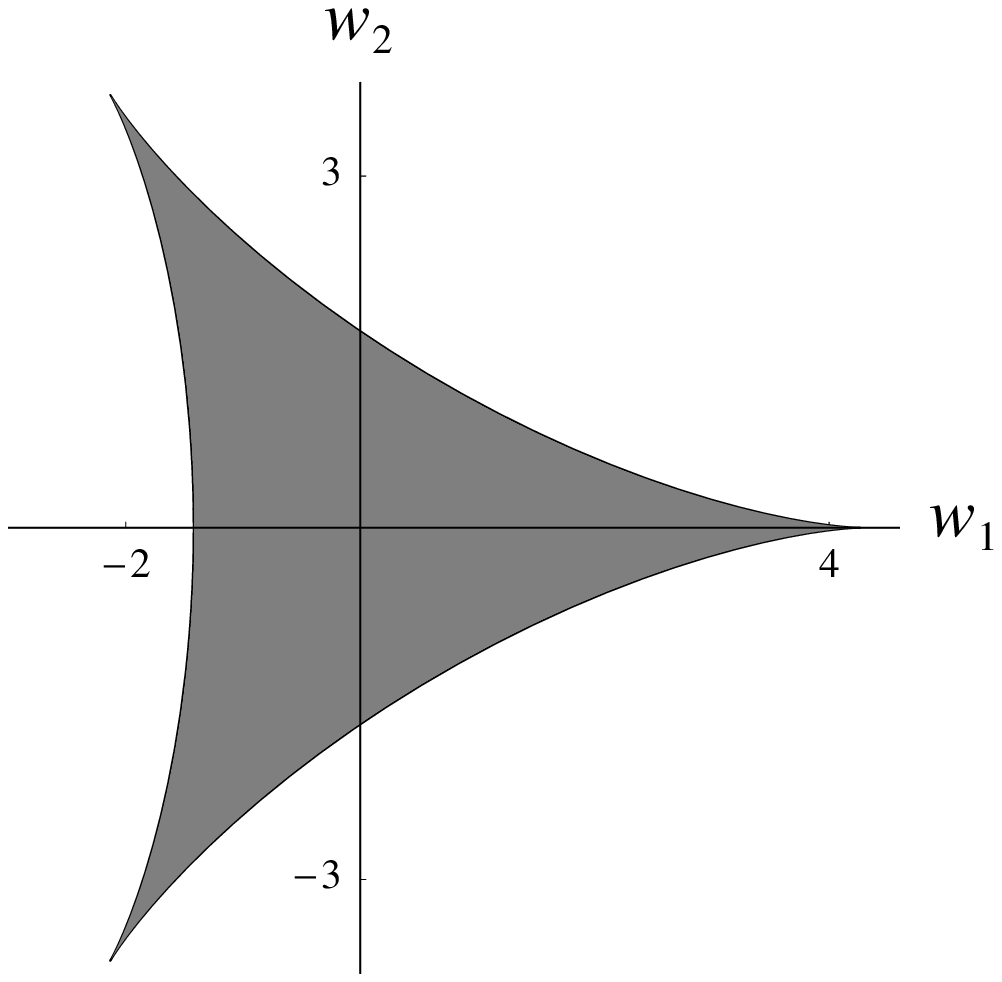}
\begin{caption}{$W_3$, the set of $(w_1, w_2)$ such that $x^3 + (w_1 +
  i w_2) x^2 + (w_1 - i w_2) x + 1$ has all roots on $\T$.} 
\label{fig:W3}
\end{caption}
\end{figure}

Given $\mathbf{a} \in \R^{N-1}$, $X_N \mathbf{a}$ is a vector in
$\C^{N-1}$.  Moreover if $\mathbf{c} = X_N \mathbf{a}$ then $c_{N-n} =
\overline{c_{n}}$ for $1 \leq n \leq N-1$, and we may associate a CR
polynomial to $\mathbf{a}$ by specifying that
\[
\mathbf{a}(x) = (x^N + 1) + \sum_{n = 1}^{N-1} c_n x^{N-n}.
\]
Let $\Delta \subset \C$ denote the open unit ball and let $\T \subset
\C$ denote the unit circle.  Equation (\ref{eq:2}) implies that there
exist $\alpha_1, \alpha_2, \ldots, \alpha_M \in \Delta$ and $\xi_1,
\xi_2, \ldots, \xi_L \in \T$ such that
\begin{equation}
\label{eq:15}
\mathbf{a}(x) = \prod_{m=1}^M (x - \alpha)\left(x -
1/\overline{\alpha}\right)\prod_{l=1}^L (x - \xi_l),
\end{equation}
where obviously $2M + L=N$.  We define $W_N$ to be the set
\[ W_N= \left\{\mathbf{w} \in \R^{N-1} : \mathbf{w}(x) \mbox{ has all
roots on } \T \right\},
\]
so that $W_N$ is in one-to-one correspondence with the set of CR
polynomials of degree $N$ with all roots on the unit circle.
Figures~\ref{fig:W3} and \ref{fig:W4} (p.\pageref{fig:W4}) show $W_3$
and $W_4$ respectively.  Elements of $W_N$ will be regarded as either
CR polynomials or as vectors in $\R^{N-1}$ as is convenient.

From Figure~\ref{fig:W3} we see that $\partial W_3$ has a distinctive
shape.  In fact, it is a 3-cusped hypocycloid. (A hypocycloid is the
curve produced by the image of a point on a circle as the circle rolls
around the inside of a larger circle.)  Similarly, the projection of
$W_4$ onto the $w_1, w_3$ plane is bounded by a 4-cusped hypocycloid.

\subsection{Statement of Results}

The geometric properties of $W_3$ and $W_4$ (figures~\ref{fig:W3} and \ref{fig:W4}) suggest
many patterns in the structure of $W_N$ in general.  
\begin{thm}
\label{thm:1} $W_N$ is homeomorphic to $B^{N-1}$, the closed $N-1$
dimensional ball.
\end{thm}
We will see that $W_N$ is circumscribed by the sphere of radius
$\sqrt{ {2N \choose N} -2}$. Moreover, the CR polynomials
$(x+\zeta_N)^N$, where $\zeta_N$ is an $N$th root of unity, correspond
to the only points in $W_N$ intersecting this sphere.  Given
$\mathbf{w} \in W_N$, there is a natural way to define a (cyclically
ordered) partition of $N$, $\mathcal{P}(\mathbf{w})$, corresponding to
the multiplicities of the cyclically ordered roots of $\mathbf{w}(x)$.
In this manner, $W_N$ has a finer structure imposed upon it.
\begin{cor}
\label{thm:2} $W_N$ has the structure of a coloured $N-1$ simplex
  where the colouring of $\mathbf{w}\in W_N$ is
  $\mathcal{P}(\mathbf{w})$.
\end{cor}
The polynomials $(x+\zeta_N)^N$ are vertices in the above
correspondence.  This colouring seems to affect the geometry of $W_N$.
For example, the edges of $W_4$ have different curvatures, depending
on their colouring, and one can show that, despite appearances, $W_4$
is not star-shaped with respect to the origin.  These geometric constraints
are also evident in the isometry group of $W_N$, which is a proper
subgroup of the group of isometries of an $N$-simplex.
\begin{thm}
\label{thm:3} The group of isometries of $W_N$ is isomorphic to $D_N$.
As such, $\Isom(W_N)$ is generated by $R$ and $C$, where for any
$\mathbf{w} \in W_N$, $R$ corresponds to the action of multiplying
each root of $\mathbf{w}(x)$ by $\zeta_N$, and $C$ corresponds to
complex conjugation of the roots of $\mathbf{w}(x)$.  That is, if
$\mathbf{w}(x) = \prod_{n=1}^N (x - \xi_n)$, then
$$
R \cdot \mathbf{w}(x) = \prod_{n=1}^N (x - \zeta_N \xi_n) \quad
\mbox{and} \quad C \cdot \mathbf{w}(x) = \prod_{n=1}^N (x -
\overline{\xi_n}).
$$
\end{thm}
We will also explicitly compute the volume of $W_N$ by using techniques from 
random matrix theory. By the {\it volume} of $W_N$ we simply mean its Lebesgue
measure.
\begin{thm} 
\label{thm:4}
The volume of $W_N$ is given by
\[ 
\vol(W_N) = \frac{2^{N-1} \pi^{(N-1)/2}}{\Gamma(\frac{N+1}{2})}.
\]
That is, the volume of $W_N$ is equal to the volume of the $N-1$
dimensional ball of radius $2$.
\end{thm}
Notice that $\vol(W_N)$ is always a rational number times an integer
power of $\pi$, since when $M$ is odd, $\Gamma(M)$ is a rational
number times $\sqrt{\pi}.$

In Section~\ref{preliminary} we will prove some preliminary results
about $W_N$.  Theorem~\ref{thm:1} and Corollary~\ref{thm:2} will be
proved in Section~\ref{topology}.  Theorem~\ref{thm:3} will be proved
in Section~\ref{geometry}, and Theorem~\ref{thm:4} will be proved in
Section~\ref{volume}.

\subsection{Motivation}

Conjugate reciprocal polynomials have appeared in the literature in
the study of random polynomials, random matrix theory and most
recently speculative number theory.  The primary interest in CR
polynomials, and in particular polynomials in $W_N$ has been in the
distribution of their roots.  The statistical behavior of the roots of
polynomials in $W_N$ in the large $N$ limit has been used
in the study of quantum chaotic dynamics \cite{MR1418808} and the
distribution of zeros on the critical line of certain $L$-functions.
In this vein, D.~Farmer, F.~Mezzadri and N.~Snaith have studied the zeros
of polynomials in $W_N$ as a model for the distribution of zeros of
$L$-functions which do not have an Euler product but nonetheless
satisfy the Riemann Hypothesis \cite{FMS05}.  The second author was 
introduced to the problem of finding the volume of $W_N$ by D. Farmer.

Likewise reciprocal polynomials have been of interest in number theory
to those studying Mahler measure \cite{smyth} and those studying 
abelian varieties over finite fields.  S.~DiPippo and E.~Howe use the
volume of the set of monic real polynomials of degree $N$ with all
roots on the unit circle to give asymptotic estimates for the number
of isogeny classes of $N$-dimensional abelian varieties over a finite
field \cite{dipippo-howe}.  A monic polynomial with real coefficients
and all roots on the unit circle is necessarily reciprocal and hence
the set of all such polynomials of degree $N$ is a set akin to $W_N$.
The set of non-monic reciprocal polynomials with Mahler measure
equal to one and degree at most $N$ has been used to give asymptotic
estimates for the number of reciprocal polynomials in $\Z[x]$ of
degree at most $N$ with Mahler measure bounded by $T$ as $T
\rightarrow \infty$ \cite{sinclair-2005}.  The volume of the set of
non-monic polynomials in $\R[x]$ of degree $N$ with all roots in the
closed unit disk has been used to give asymptotic estimates for the
number of polynomials in $\Z[x]$ of degree at most $N$ with Mahler
measure bounded by $T$ as $T \rightarrow \infty$ \cite{chern-vaaler}.

This manuscript will be largely concerned with the global properties
of the coefficient vectors of CR polynomials in $W_N$.  As we shall
see, the CR condition together with the condition that all roots lie on
the unit circle constrains the geometry and topology of $W_N$.  Here
we study the geometry, topology and volume of $W_N$ not for any
particular application, but for their own sake.

\subsection{Preliminary results about $W_N$}\label{preliminary}

The following proposition gives a useful characterization of elements
of $W_N$ based on their roots.

\begin{prop}
\label{CR characterization prop} The vector $\mathbf{w}$ is in $W_N$
if and only if
$$
\mathbf{w}(x) = \prod_{n=1}^N (x - \xi_n),
$$
where $\xi_1, \xi_2, \ldots, \xi_N$ are elements of $\T$ satisfying $\xi_1 \xi_2
\cdots \xi_N = (-1)^N$.
\end{prop}
\begin{proof} Suppose $\mathbf{w}$ is as in the statement of the
proposition.  Then,
\begin{eqnarray*} x^N \overline{\mathbf{w}(1/\overline{x})} &=& x^N
\prod_{n=1}^N \left(\frac{1 - \overline{\xi_n} x}{x} \right) = (-1)^N
\prod_{n=1}^N (\xi_n^{-1} x - 1) = \prod_{n=1}^N \xi_n (\xi_n^{-1} x -
1).
\end{eqnarray*} It follows that $\mathbf{w}(x) = x^N
\overline{\mathbf{w}(1/\overline{x})}$, and thus $\mathbf{w} \in W_N$.
The converse is obvious since every element of $W_N$ is a polynomial
with all roots on the unit circle and constant coefficient 1.
\end{proof}

In order to exploit both the symmetry of the coefficients of
polynomials in $W_N$ and the symmetry of the roots we introduce the
map $E_N: \C^N \rightarrow \C^N$ given by 
$E_N(\boldsymbol{\alpha}) = \mathbf{b}$ where
$$
x^N + \sum_{n=1}^N b_{N-n} x^n = \prod_{n=1}^N (x - \alpha_n).
$$
That is, the $n$-th coordinate function of $E_N(\boldsymbol{\alpha})$
is given by $(-1)^n e_n(\alpha_1, \alpha_2, \ldots, \alpha_N)$.
Clearly, since the coefficients of a monic polynomial are independent
of the ordering of the roots, $E_N$ induces a map (which we also
call $E_N$) from $\C^N / S_N$ to $\C^N$. (Where $S_N$ is the symmetric
group on $N$ letters, and $\C^N / S_N$ is the orbit space of $\C^N$
under the action of $S_N$ on the coordinates of $\bs{\alpha}$.)  The
torus $\T^N$ sits in $\C^N$, and thus $E_N$ gives a correspondence
between the $\T^N / S_N$ and the set of polynomials of degree $N$ in
$\C[x]$ with all roots on the unit circle.  We define $\Omega_N$ to be
the subset of $\T^N / S_N$ given by
$$
\Omega_N = \left\{(\xi_1, \xi_2, \ldots, \xi_{N}) : \xi_1 \xi_2 \cdots
\xi_N = (-1)^N \right\}.
$$
If $\mathbf{w} = E_N(\boldsymbol{\xi})$ then $\boldsymbol{\xi}$ will
be referred to as a {\it root vector} of $\mathbf{w}(x)$.  Clearly, by
Proposition \ref{CR characterization prop}, the map $E_N$ induces a
homeomorphism between $\Omega_N$ and $W_N$.

We now turn to the structure of $W_N$ viewed as a subset of
$\R^{N-1}$.   

\begin{prop} 
\label{prop:1}
$W_N$ is a closed path connected set with positive
  volume.  Moreover, the boundary of $W_N$ is given by 
$$
\partial W_N = \{ \mathbf{w} \in W_N : \disc(\mathbf{w}) = 0 \},
$$
where $\disc(\mathbf{w})$ is the discriminant of $\mathbf{w}(x)$.
\end{prop}
\begin{proof} 
Suppose $\mathbf{a} \in \R^{N-1}$, and $\mathbf{a} \not\in W_N$.
Then, by (\ref{eq:15}), there exists $\alpha \in \Delta$ such that
$(x - \alpha)(x - 1/\overline{\alpha})$ is a factor of
$\mathbf{a}(x)$.  Since the root vector of a polynomial is a
continuous function of the roots, there must exist some open
neighborhood $U$ of $\mathbf{a}$ such that $U \cap W_N = \emptyset$.
It follows that $W_N$ is closed.

Now suppose $\mathbf{w} \in W_N$ has a double root, that is there exists
$\xi_1, \xi_2, \ldots, \xi_{N-1} \in \T$ such that
$\mathbf{w}(x) = (x - \xi_1)^2 (x - \xi_2) \cdots (x - \xi_{N-1})$.
If $U$ is an open neighborhood of $\mathbf{w}$ then, since the
coefficients of a 
polynomial are continuous functions of the roots, there exists
$\epsilon > 0$ such that $(x - \epsilon \xi_1)(x -
\epsilon^{-1} \xi_1)(x - \xi_2) \cdots (x - \xi_{N-1})$ is a CR
polynomial in $U \setminus W_N$.  That is, $\mathbf{w} \in \partial
W_N$.  It follows that $\partial W_N$ consists of polynomials with
$\disc \mathbf{w} = 0$.

To see that $W_N$ has positive volume, let $\mathbf{v}$ be in
$\overcirc{W_N},$ the interior of $W_N$. Then $\disc(\mathbf{v}) \neq
0$ and since the discriminant of a polynomial is a continuous function
of the coefficients, there must exist an open neighborhood $U$ of
$\mathbf{v}$ such that all polynomials in $U$ have non-zero
discriminant.  Thus $U \subset \overcirc{W_N}$ and $W_N$ has positive
volume.

To see that $W_N$ is path connected, let $\mathbf{w}_1$ and
$\mathbf{w}_2$ be two points in $W_N$. Marking $N-1$ roots in each
$\mathbf{w}_1$ and $\mathbf{w}_2$ it is clear that there is a
continuous map between the $N-1$ dimensional root vectors where all
roots are in $\T$. By
stipulating that the final roots are assigned to satisfy the CR
condition, we see that this extends to a continuous map with image in
$W_N$ from $\mathbf{w}_1$ to $\mathbf{w}_2$.
\end{proof}

We now turn to the geometry of $W_N$.
Our first result in this direction is identifying those points in $W_N$
which are farthest from the origin, which corresponds to the CR polynomial $x^N+1$.
\begin{prop}
\label{distance prop} If $\mathbf{w} \in W_N$, then
$$\| \mathbf{w} \|^2 \leq {2 N \choose N} - 2.$$
Moreover, there is equality in this equality if and only if
$\mathbf{w}(x) = (x + \zeta_N)^N$, where $\zeta_N$ is an $N$th root of
unity.
\end{prop}
\begin{proof}
Suppose $\mathbf{w} \in W_N$ and that $\mathbf{w}(x) = \prod_{n=1}^N
(x - \xi_n),$ for $\xi_1, \xi_2, \ldots, \xi_N \in \T$.  Let
$\mathbf{b} = X_N \mathbf{w}$ and notice that
\[
| b_n |^2 =\left| \sum_{\mathbf{i} \in \mathcal{I}(N,n) } \left\{
\prod_{m = 1}^n \xi_{i_m} \right\} \right|^2,
\]
where $\mathcal{I}(N,n) = \left\{ \mathbf{i} \in \Z^n : 1 \leq i_1 <
i_2 < \cdots < i_n \leq N \right\}$.  Clearly then, 
\[
| b_n |^2 \leq {N \choose n}^2
\quad \mbox{with equality when} \quad 
\prod_{m=1}^n \xi_{i_m} = \prod_{m=1}^n \xi_{j_m} 
\]
for every choice of $\mathbf{i}, \mathbf{j} \in \mathcal{I}(N,n)$.
That is, there is equality exactly when $\xi_1 =
\xi_2 = \cdots = \xi_N$.  In this case, Proposition \ref{CR
characterization prop} implies that $\xi_1^N = (-1)^N$, which is only
satisfied when $\xi_1 = \xi_2 = \cdots = \xi_N = \zeta_N$ for some
$N$th root of unity.

It follows that
$$
\|\mathbf{w}\|^2 = \| \mathbf{b} \|^2 \leq
\sum_{n=1}^{N-1} {N \choose n}^2 = {2 N \choose
N} - 2,
$$ 
where the last equality comes from the well known formula for the sum
of the squares of binomial coefficients.  By our previous remarks,
equality is attained in the inequality only when $\mathbf{w}(x)$ is a
polynomial of the form $(x + \zeta_N)^N$.
\end{proof}

This proposition gives us a hint of the geometric structure of $W_N$.
Let $\zeta_N = e^{2 \pi i/N}$, and let $\mathbf{v}_1, \mathbf{v}_2,
\ldots, \mathbf{v}_N \in \R^{N-1}$ be determined by setting
$\mathbf{v}_n(x) = (x + \zeta_N^n)^N$.  In particular,
$\mathbf{v}_N(x) = (x + 1)^N$.  For reasons which will become clear
we will call $\mathbf{v}_1, \mathbf{v}_2, \ldots, \mathbf{v}_N$ the
{\it vertices} of $W_N$.  As we will see, every isometry of $W_N$ must
fix the origin.  This together with Proposition \ref{distance prop}
implies that the set of isometries of $W_N$ must permute its vertices.
We will also prove that the vertices span 
$\mathbb{R}^{N-1}$.  It follows that  if two isometries $T_1$ and $T_2$
induce the same permutation of the vertices, the isometry
$T_1T_2^{-1}$ fixes every vertex.  Therefore, since $T_1T_2^{-1}$
extends to an isometry that fixes a spanning set for $R^{N-1}$, it is
the identity.  Therefore, the group of isometries of $W_N$ is
isomorphic to a subgroup of $S_N$, as each isometry is uniquely
determined by the permutation it induces on the set of vertices.  In
fact, Theorem~\ref{thm:3} shows that this group is isomorphic to $D_N$, the group of
isometries of a regular $N$--gon.

At present, for $R$ and $C$ defined as in Theorem~\ref{thm:1} we will demonstrate
\begin{prop}\label{RC}
$R$ and $C$ are isometries of $W_N$.
\end{prop}
We  defer the proof that they generate the complete group of
isometries until Section~\ref{sec:proof-theor-refthm:4}.  

\begin{proof}
We use $\|
\cdot \|$ to denote the usual 2--norm on both $\R^{N-1}$ and
$\C^{N-1}$.  Given $\mathbf{w}_1,\mathbf{w}_2 \in W_N$, it is easily
verified that $\|X_N \mathbf{w}_1\| = \| \mathbf{w}_1 \|$.  Thus,
since $X_N \mathbf{w}_1 - X_N \mathbf{w}_2$ is the coefficient vector
of a polynomial, Parseval's formula yields
\begin{equation}
\label{eq:12}
\| \mathbf{w}_1 - \mathbf{w}_2\|^2 =  \| X_N \mathbf{w}_1 - X_N
\mathbf{w}_2\|^2 = \frac{1}{2\pi} \int_N^{2 \pi} \left|
    \mathbf{w}_1(e^{i \theta}) - \mathbf{w}_2(e^{i\theta}) \right|^2 \,
  d\theta.
\end{equation}
Notice that $R \cdot
\mathbf{w}_1(x) = \mathbf{w}_1(\zeta_N^{-1} x)$, and hence
\begin{eqnarray*} \|R \cdot \mathbf{w}_1 - R \cdot \mathbf{w}_2\|^2 &=&
 \frac{1}{2 \pi} \int_N^{2 \pi} \left| R \cdot
\mathbf{w}_1(e^{i \theta}) - R \cdot \mathbf{w}_2(e^{i \theta}) \right|^2
\, d\theta \\ &=&\frac{1}{2 \pi}
\int_N^{2 \pi} \left| \mathbf{w}_1(\zeta_N^{-1} e^{i \theta}) -
\mathbf{w}_2(\zeta_N^{-1} e^{i\theta}) \right|^2 \, d\theta =
\| \mathbf{w}_1 - \mathbf{w}_2 \|^2 \\
\end{eqnarray*} 
where the last equation follows from an easy change of variables.
Similarly, $C \cdot \mathbf{w}_1(x) =
\overline{\mathbf{w}_1(\overline{x})}$ and thus,
\begin{eqnarray*} \|C \cdot \mathbf{w}_1 - C \cdot \mathbf{w}_2\|^2 &=&
 \frac{1}{2 \pi} \int_N^{2 \pi} \left| C \cdot
\mathbf{w}_1(e^{i \theta}) - C \cdot \mathbf{w}_2(e^{i \theta}) \right|^2
\, d\theta \\ &=& \frac{1}{2 \pi}
\int_N^{2 \pi} \left| \mathbf{w}_1(e^{-i \theta}) - \mathbf{w}_2(e^{-i
\theta}) \right|^2 \, d\theta = \| \mathbf{w}_1 - \mathbf{w}_2
\|^2. 
\end{eqnarray*} 
It follows that $R$ and $C$ are isometries of $W_N$.
\end{proof}

\subsection{The Colouring of $W_4$}\label{W4}

As mentioned previously, every $\mathbf{w} \in W_N$ uniquely
determines a partition $\mathcal{P}(\mathbf{w})$ of $N$ up to cyclic
ordering, corresponding to the multiplicities of the cyclically
ordered roots of $\mathbf{w}(x)$. That is, we may decompose $W_N$ into
regions, {\it faces} if you will, determined by partitions of the
integer $N$.

\begin{figure}[h!]
\centering
\includegraphics[scale=.5]{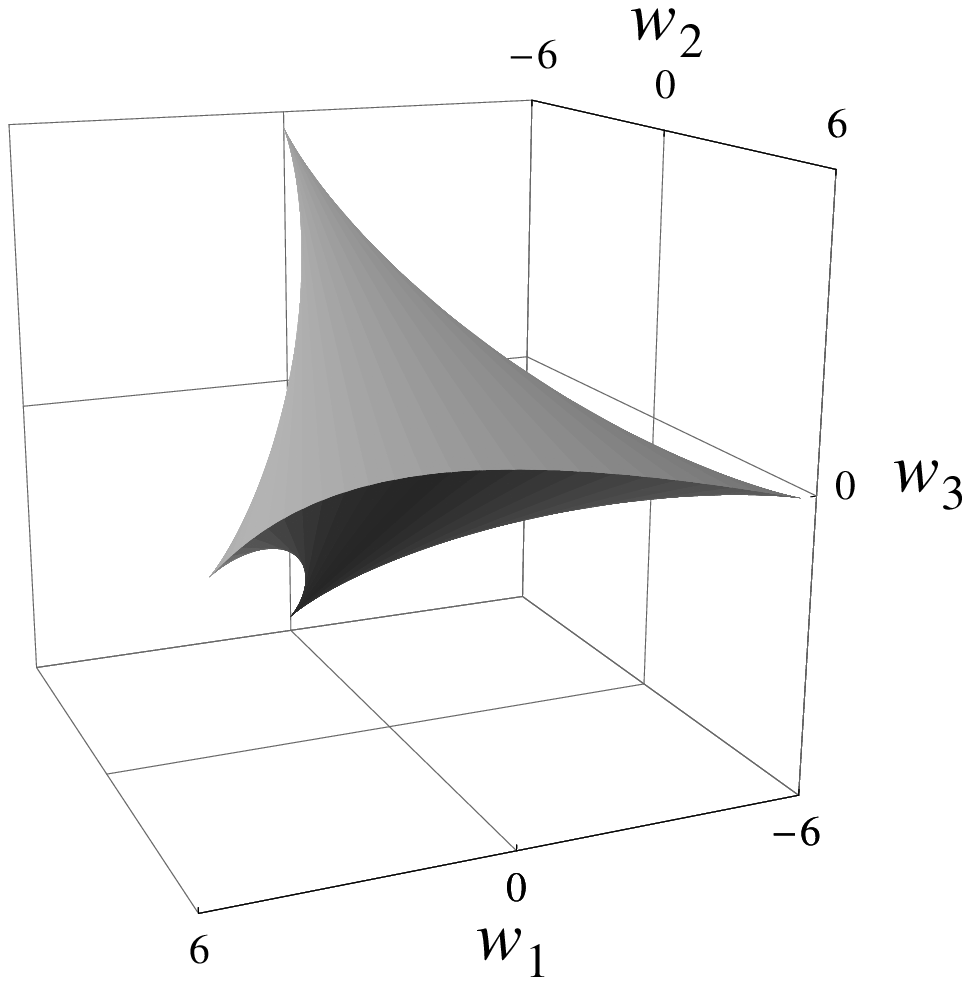}
\includegraphics[scale=.5]{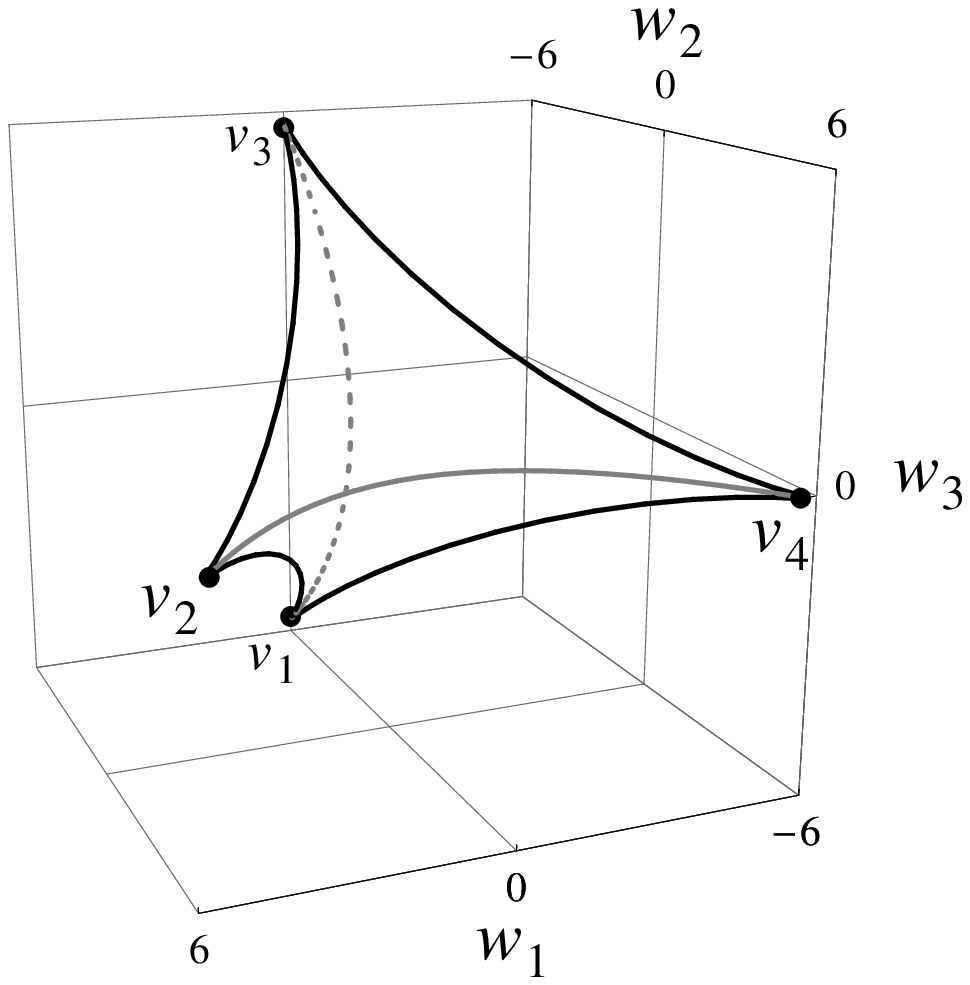}
\begin{caption}{The faces of $W_4$} 
\label{fig:W4}
\end{caption}
\end{figure}
For example, Figure~\ref{fig:W4} demonstrates the decomposition of
$W_4$ into a coloured simplex. The interior is associated to the partition 
$(1,1,1,1)$.  The four vertices, given by $\mathbf{v}_k(x) = (x-i^k)^4$ for 
$1 \leq k \leq 4$, each correspond to the partition $(4)$.  Likewise, 
there are four faces, each associated to the partition 
$(2,1,1)$.  In contrast, there are  two different partitions of $4$ of
length two: $(2,2)$ and $(3,1)$.  Each of these partitions corresponds
to faces of codimension $2$, that is {\it edges} of $W_4$.  There are
six edges, four of type $(3,1)$ (black curves on the right in
Figure~\ref{fig:W4}) and two of type $(2,2)$ (gray curves).  One of
the $(2,2)$ edges joins $\mathbf{v}_2$ with $\mathbf{v}_4$ and the
other joins $\mathbf{v}_1$ with $\mathbf{v}_3$.

To illustrate the difference between edges
corresponding to the partition $(2,2)$ and edges corresponding to
$(3,1)$, consider a path starting at $\mathbf{v}_4$ and traversing an
edge to another vertex.  The polynomial corresponding to
$\mathbf{v}_4$ has a root of multiplicity four at $x=1$.  The edge
associated to the partition $(2,2)$ consists of polynomials with two
double roots.  If we imagine one of these double roots starting at
$x=1$ and traversing the unit circle in the counterclockwise
direction, then the other double root must traverse the circle in the
clockwise direction.  Moreover, these double roots must traverse the
circle at the same rate so that each intermediate polynomial on the
$(2,2)$ edge has constant coefficient equal to 1.  A moment's thought
reveals that the only other vertex which can be reached from
$\mathbf{v}_4$ via a $(2,2)$ edge corresponds to the polynomial with a
root of multiplicity four at $x=-1$, i.e.~$\mathbf{v}_2$.  The other
$(2,2)$ edge connects $\mathbf{v}_1$ with $\mathbf{v}_3$.  On the
other hand, if we look at an edge starting at $\mathbf{v}_4$ formed by
dividing the root of multiplicity four into a root of multiplicity
three traversing the circle in the counterclockwise direction and a
single root traversing the circle in the clockwise direction, the
single root must traverse the circle at a rate three times that of the
triple root so that the intermediate polynomials have constant
coefficient equal to $1$.  We see that when the triple root has
reached $x=i$ then the single root too is at $x=i$ and thus a
$(3,1)$ edge connects $\mathbf{v}_4$ with $\mathbf{v}_1$.  Similarly,
when the triple root reaches $x=-1$ then so too has the single root,
so that a $(3,1)$ edge connects $\mathbf{v}_1$ with $\mathbf{v}_2$.
Continuing in this manner, we see that any two vertices can be
connected via an edge consisting of polynomials with a triple root and
a single root.

The group generated by $R$ acts transitively on both the set of
$(2,2)$ edges and the set of $(3,1)$ edges.  This, coupled with the
action of $C,$ demonstrates that there is an isometry sending any edge
to itself which interchanges the vertex endpoints.  Thus the curvature
of each edge is symmetric.  No isometry of $W_N$ carries a $(2,2)$
edge to a $(3,1)$ edge.  Indeed, the curvature of the $(2,2)$ edges
differs from the curvature of the $(3,1)$ edges.  This reflects the
fact that $\Isom(W_4) \cong D_4$, not the full symmetric group.

\section{The Topology of $W_N$}\label{topology}

\subsection{The Proof of Theorem~\ref{thm:1}} 
Recall that $W_N$ is homeomorphic to $\Omega_N $ where $\Omega_N$ is
the subset of $\mathbb{T}^N/S_N$ given by
$$ \Omega_N=\left\{(\xi_1, \xi_2, \ldots, \xi_{N}) : \xi_1 \xi_2 \cdots
	\xi_N = (-1)^N \right\}. $$ 
Let $c=0$ if $N$ is even and $c=1/2$ if $N$ is odd.	
	After reparametrization by
$\xi_n=e^{2\pi i \vartheta_n}$ we have
\begin{equation}
\label{eq:6}
W_N \cong \left\{(\vartheta_1, \vartheta_2, \ldots, \vartheta_{N}) : 
          \left\{ \sum_{n=1}^N \vartheta_n \right\} = c \right\}
\end{equation}
where $\vartheta_n \in \mathbb{R}/\mathbb{Z}$, for $1\leq n \leq N$,
and in this context $\{x\}$ denotes the fractional part of $x$.  We
will denote the set on the right hand side of (\ref{eq:6}) by $\Theta_N$,
and the interior of this set $\overcirc{\Theta}_N$.  Throughout this
discussion, we will continue to refer to the $\vartheta_n$ as
roots, and we will denote the reparameterized torus $\R /\Z$ as
$\T_{+}$. (The subscript reflects the fact that in this context
we are working with the additive torus). 

First, we will prove that the interior of $W_N$ is homeomorphic to an
open $N-1$ ball.  Define $\Phi: \overcirc{\Theta}_N\rightarrow
\T_{+}^{N-1}/S_{N-1}$ as follows.  Fix a basepoint
$\boldsymbol{\psi}=(\psi_1, \dots, \psi_N) \in \overcirc{\Theta}_N$.
Since $\boldsymbol{\psi}$ corresponds to a polynomial in the interior
of $W_N$, by Proposition~\ref{prop:1}, all coordinates of
$\boldsymbol{\psi}$ are unique.  Fix some $n\in \{1, \dots, N\}$.  We
define $\Phi(\boldsymbol{\psi})=(\psi_1, \dots, \hat{\psi_n}, \dots
\psi_N)$, where $\hat{\psi_n}$ means that $\psi_n$ is omitted.  We can
extend $\Phi$ in a neighborhood of $\boldsymbol{\psi}$ since the roots
vary continuously.  Specifically, for
$\boldsymbol{\vartheta}=(\vartheta_1, \dots, \vartheta_N),$
$\Phi(\boldsymbol{\vartheta})=(\vartheta_1, \dots, \hat{\vartheta_l},
\dots, \vartheta_N)$ where $\vartheta_l$ is the closest root of to
$\psi_n$ on $\T_{+}$.  In fact, for any $\boldsymbol{\vartheta} \in
\overcirc{\Theta}_N$, given a path in $\overcirc{\Theta}_N$ from
$\boldsymbol{\psi}$ to $\boldsymbol{\vartheta}$ we can extend $\Phi$
to $\boldsymbol{\vartheta}$ due to the continuity of the roots and the
fact that there are no multiple roots in the interior of $W_N$.  

To define $\Phi$ on all of $\overcirc{\Theta}_N$ it is enough to show that
$\Phi$ does not depend on the choice of path.  It suffices to show
that if $\gamma :[0,1] \rightarrow \overcirc{\Theta}_N$ is a
continuous loop based at $\boldsymbol{\vartheta}$, then $\Phi \circ
\gamma(0)=\Phi \circ \gamma(1)$.  Since there are no double roots in
$\gamma(t)$, the continuous orbit of the roots of
$\boldsymbol{\vartheta}$ on $\T_{+}$ by $\gamma$ is homotopic to a
rotation.  As such, there are continuous functions
$\tau_n:[0,1] \rightarrow \R$ such that
$\gamma(t)=(\vartheta_1+\tau_1(t),\dots,\vartheta_N+\tau_N(t))$.  The
conjugate reciprocal condition translates to the condition
\[
\left\{ \sum_{n=1}^N (\vartheta_n+\tau_n(t))\right\} = c 
\]
where $c$ is as before.  As 
\[
\left\{ \sum_{n=1}^N \vartheta_n \right\} = c 
\qquad
\mbox{we conclude that}
\qquad
\left\{ \sum_{n=1}^N \tau_n(t)\right\} = 0.
\]

Fixing representatives $\vartheta_1, \dots , \vartheta_N$ in $[0,1)$,
we define $\sigma: [0,1] \rightarrow \R$ by
\[
\sigma(t)=\sum_{n=1}^N (\vartheta_n+\tau_n(t)).
\]
The continuity of $\sigma$ follows from the continuity of the
$\tau_n$.  Since $\gamma$ is homotopic to a rotation, there is a
unique $1 \leq k'\leq N$ such that $\gamma(\vartheta_k) =
\vartheta_{k+k' \bmod N}$ for $1 \leq k \leq N$. It follows that
$\sigma(1) = \sigma(0)+{|k'|}$.  As $\{ \sigma(t)\}= 0$ for all $t\in
[0,1]$, this contradicts the continuity of $\sigma$ unless $k'\equiv 0
\bmod{N}$.  Hence $\Phi$ is well-defined on $\overcirc{\Theta}_N$ and
can extended to all of $\Theta_N$.  Moreover, $\Phi$ is injective,
since if $\Phi(\boldsymbol{\vartheta})=\Phi(\boldsymbol{\vartheta}')$
then $\boldsymbol{\vartheta}$ and $\boldsymbol{\vartheta}'$ share
$N-1$ roots and by the conjugate reciprocal condition,
$\boldsymbol{\vartheta}=\boldsymbol{\vartheta}'$.  The conjugate
reciprocal condition allows us to define an inverse of $\Phi$ and, as
$\Phi$ is continuous, we conclude that $\overcirc{\Theta}_N$ is
homeomorphic onto a subset of $\T_{+}^{N-1}/S_{N-1}$.

We now shift our attention to $\partial \Theta_N$; first we introduce
some notation.  Given positive integers $M, m$ and $n$ such that $1
\leq n < m \leq M$, define
\[
L_m^M=\Bigg\{ (\vartheta_1, \dots , \vartheta_{M})\in \T_{+}^{M}
			: \Bigg\{ 2\vartheta_n+\sum_{l=1, l \neq
                            m}^{M} \vartheta_l \Bigg\}= c \Bigg\}
\]
and 
\[
K_{n,m}^M=\left\{  (\vartheta_1, \dots , \vartheta_{M})\in \T_{+}^M
			: \vartheta_n=\vartheta_m\right\}.
\]
We then set
\[
L^M=\bigcup_{m=1}^M L_m^M, \qquad K^M=\bigcup_{1\leq
  n<m\leq M} K_{n,m}^M \qquad \mbox{and} \qquad J^M = L^M \cup K^M.
\]
Clearly $L^M, K^M$ and $J^M$ are stabilized by $S_M$, and thus (for
instance) $(\T^{M}_+ \setminus J^{M})/S_{M} = (\T^{M}_+ / S_{M})
\setminus (J^{M} / S_{M})$.  Returning to $\Phi$, if a point in
$J^{M-1}/S_{N-1}$ is in the image of $\Phi$, then its preimage
necessarily has multiplicity at least two. Therefore
$\Phi(\overcirc{\Theta}_N)$ is a subset of $(\T_{+}^{N-1} \setminus
J^{N-1})/S_{N-1}$ and $\Phi$ maps the boundary of $\Theta_N$ into
$J^{N-1}/S^{N-1}$.

\begin{figure}[h]
\centering
\includegraphics[scale=.8]{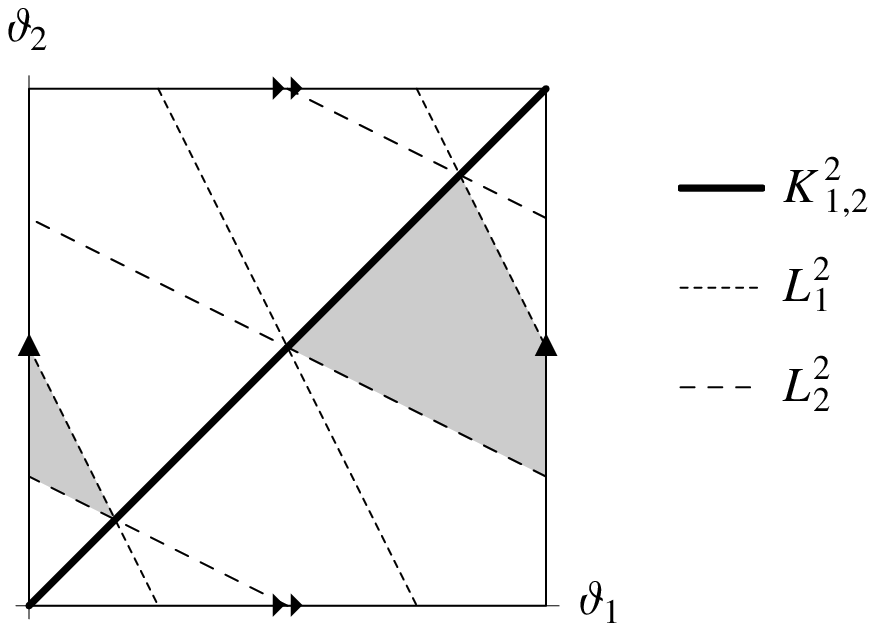}
\begin{caption}{
An open 1-ball in $\T_{+}^{2}\setminus J^{2}$
}
\label{fig:torusfun}
\end{caption}
\end{figure}

For any positive $M$, the removal of a hyperplane from $\T_{+}^M$
corresponds to the removal of a generator in the fundamental group.
Since the hyperplanes $L_{m,n}^M$ are liearly independent, the
fundamental group of $\T_+^M \setminus L^M$ is trivial. 
We conclude that the set $\mathbb{T}^M
\setminus L^M$ is homeomorphic to a disjoint union of open $M$-balls.
Therefore, the deletion from $\T_{+}^M \setminus L^M$ of the set of
hyperplanes $K^M$ is also homeomorphic to a disjoint union of open
$M$-balls. Figure~\ref{fig:torusfun} shows one of these open balls in
$\T_{+}^2 \setminus J^2$. The action by $S_M$ on
$\T_{+}^M$ is given by reflections through the $K_{m,n}^M$. As such,
since $K^M$ and $L^M$ are stabilized by this action, any reflection
through $K_{m,n}^M$ must map a ball onto another ball. The image of
$\overcirc{\Theta}_N$ under $\Phi$ must therefore be one of these
balls as $\partial \Theta_N$ maps into $J^{N-1}$, and we
conclude that $\Phi(\Theta_N)$ is homeomorphic to a closed
$(N-1)$-ball.  Since $\Phi$ is a homeomorphism, we conclude that $\Theta_N$ and therefore $W_N$ is homeomorphic to a closed $(N-1)$-ball.

\subsection{The Proof of Corollary~\ref{thm:2}}

Let $\mathcal{P}$ be a partition of $N$, {\it i.e.}~$\mathcal{P} =
(n_1 \dots, n_M)$ is a vector of positive integers such that $n_1 +
n_2 + \cdots + n_M = N$.  Let $|\mathcal{P}| = M$ denote the length of
the partition.  By associating the cyclically ordered roots of
$\mathbf{w}(x)$ to their multiplicities, $\mathbf{w} \in W_N$
determines a partition $\mathcal{P}(\mathbf{w})$ of $N$ which is
well-defined up to cyclic ordering.  (We will consider all partitions
only up to cyclic ordering.)  If $\mathcal{P} =(n_1, \dots, n_M)$ and
$\mathcal{P}'=(n_1+n_2, \dots, n_M)$ then we say that $\mathcal{P}'$
is obtained from $\mathcal{P}$ by reduction. Notice that reduction
gives a partial ordering on the set of cyclically ordered partitions,
and if the partition $\mathcal{P}''$ is obtained from $\mathcal{P}$ by
a series of reductions we will write $\mathcal{P}'' \preceq
\mathcal{P}$.  We call each partition of $N$ a colouring, and colour
$W_N$ according to the partition type of each $\mathbf{w} \in W_N$.

If $\mathbf{w} \in \partial W_N$ and $\mathcal{P}(\mathbf{w}) = (n_1,
n_2, \ldots, n_M)$ then there exist distinct $\vartheta_1,
\vartheta_2, \ldots, \vartheta_M \in \T_+$ such that $\mathbf{w}(x) =
(x - e^{2 \pi i \vartheta_1})^{n_1} (x - e^{2 \pi i
  \vartheta_2})^{n_2} \cdots (x - e^{2 \pi i \vartheta_M})^{n_M}$. 
Letting the $\vartheta_m$ vary, the continuity of the elementary
symmetric functions implies that 
\[
\big\{ (\vartheta_1, \vartheta_2,
\ldots, \vartheta_M) :   \{n_1 \vartheta_1 + n_2 \vartheta_2 + \cdots
+ n_M \vartheta_M \} = c \big\}
\]
locally parameterizes an $M-1$ dimensional ball in $\partial W_N$
containing $\mathbf{w}$.  Let $F(\mathbf{w})$ be the maximal connected
subset of $\partial W_N$ containing $\mathbf{w}$ such that if
$\mathbf{x} \in F$ then $\mathcal{P}(\mathbf{x}) \preceq
\mathcal{P}(\mathbf{w})$.  We fix $1 \leq n \leq M$ and define
$\Phi(\boldsymbol{\vartheta}) = (\vartheta_1, \ldots,
\hat{\vartheta_n}, \ldots, \vartheta_M)$.  As in the proof of
Theorem~\ref{thm:1}, $\Phi$ has a unique extension to $F$ and we
conclude that $F$ is homeomorphic to a closed $M-1$ dimensional ball.

The proof of Theorem~\ref{thm:1} shows that $\Phi(\Theta_N)$ is
homeomorphic to a simply connected subset of $\mathbb{T}_+^{N-1}$
bounded by hyperplanes.  This is, in turn, isometric to a subset, $U$
of $\mathbb{R}^{N-1}$ bounded by hyperplanes, the boundary of which
corresponds to $(\Phi(\Theta_N) \cap J^{N-1})/S_{N-1}$.  It follows
that $U$ has the geometric structure of a polytope, and hence each $u
\in U$ lies in the interior of a face of $U$.  Moreover, the points in
the boundary of $U$ correspond to $\mathbf{w}\in W_N$ with
$\disc(\mathbf{w}) = 0$.  If $u\in U$ is in the intersection of $M$ of
these hyperplanes (i.e. $u$ is on a codimension $M$ face), then $u$
corresponds with a $\mathbf{w}\in W_N$ with
$|\mathcal{P}(\mathbf{w})|=N-M$.  In fact, this face is the image of
$F(\mathbf{w})$ as defined above.  Under this correspondence, we see
that there are $N$ vertices of the polytope, confirming that $U$ has
the structure of a simplex. It follows that $W_N$ has the structure of
a coloured simplex. Moreover, given any $\mathbf{w}\in W_N$ in the
interior of a face $F$, and any $\mathbf{u}\in \partial F$,
$\mathcal{P}(\mathbf{u}) \preceq \mathcal{P}(\mathbf{w})$.

\section{The Geometry of $W_N$}\label{geometry}

Before proving Theorem~\ref{thm:3} we need some results about the
vertices of $W_N$. 

\subsection{The Vertices of $W_N$}

\begin{lemma}
\label{lemma:1} 
The vertices of $W_N$ span $\R^{N-1}$.
\end{lemma}
\begin{proof} 
\sloppypar{
Let $U$ be the $N-1 \times N-1$ matrix whose rows are given by
$\mathbf{v}_1, \mathbf{v}_2, \ldots, \mathbf{v}_{N-1}$.  From the
definition of $\mathbf{v}_n$ it is seen that the $m,n$ entry of $U$
is given by
\[ 
U[m,n] = \left\{
  \begin{array}{ll} {\displaystyle \sqrt{2} {N \choose n} \cos\left(\frac{2\pi
          n m}{N}\right)} & \mbox{if} \quad 1 \leq n < N/2 \\  \\
    {\displaystyle {N \choose N/2} (-1)^{m-1}} & \mbox{if} \quad n =
    N/2 \\  \\ 
    {\displaystyle \sqrt{2} {N \choose n} \sin\left(\frac{2\pi (N-n)
          m}{N}\right)} & \mbox{if} \quad  N/2 < n < N.
\end{array} \right.
\] 
Hence,
\begin{equation}
\label{eq:3}
\det U = 
\left(\prod_{n=1}^{\left\lfloor N/2 \right\rfloor} {N \choose
n}\right) \left( \prod_{n=1}^{\lfloor (N-1)/2 \rfloor} {N \choose n}
\right) \det U',
\end{equation}
where $U'$ is the $N -1 \times N-1$ matrix whose $m,n$ entry is
given by
\[
U'[m,n] = \left\{
  \begin{array}{ll} \frac{\sqrt{2}}{2} \left(\e(\frac{n m}{N}) +
      \e(\frac{(N-n) m}{N})\right) & \mbox{if} \quad
    1 \leq n < N/2 \\
    \\ \e(\frac{m}{2}) & \mbox{if} n = N/2 \\
    \\ \frac{\sqrt{2}}{2} \left(i \e(\frac{n m}{N}) -i \e(\frac{(N-n)
        m}{N})\right) & \mbox{if} \quad N/2 < n < N,
\end{array} \right.
\]
and $\e(t) = e^{2 \pi i t}$.  That is, if we define the
$N-1 \times N-1$ matrix $V$ by $V[m,n] = \e(m n/N)$ then 
\begin{equation}
\label{eq:1}
U' = X_N^{\ast} V,
\end{equation}
where $X_N^{\ast}$ is the conjugate transpose of $X_N$ as defined in
(\ref{eq:8}).  This is convenient, since $|\det X_N| = 1$ and $|\det
V| = |\det V'|$ where $V'$ is the $N-1 \times N-1$ matrix given by
\[
V'[m,n] = e^{\frac{2 \pi i m n}{N} - 1}.
\]
That is $V'$ is the Vandermonde matrix formed from the complex numbers
$e(1/N)$, $e(2/N)$, $\ldots$, $e((N-1)/N)$.  The
well known relationship between Vandermonde determinants and
discriminants implies that $|\det V'| = |\disc((
x^N-1)/(x-1))|^{1/2}$.  This together with (\ref{eq:3}) and
(\ref{eq:1}) yield 
}
\[
|\det U| = 
\left(\prod_{n=1}^{\left\lfloor N/2 \right\rfloor} {N \choose
n}\right) \left( \prod_{n=1}^{\lfloor (N-1)/2 \rfloor} {N \choose n}
\right) \sqrt{\disc\left(\frac{x^N - 1}{x-1} \right)}. 
\]
Since the discriminant of $(x^N-1)/(x-1)$ is nonzero we conclude
$\det U \neq 0$ and that the vertices of $W_N$ span $\R^{N-1}$.
\end{proof}

\begin{lemma}
\label{lemma:3}
$\mathbf{v}_1 + \mathbf{v}_2 + \cdots + \mathbf{v}_N = \boldsymbol{0}.$
\end{lemma}
\begin{proof}
We use the fact that if $\zeta \neq 1$ is any $N$-th root of unity
then $\zeta + \zeta^2 + \cdots + \zeta^N = (\zeta^N - 1)/(\zeta - 1) =
0$.   If $0  < M < N$ then the coefficient of $x^{N-M}$ in
$\mathbf{v}_n(x)$ is given by ${N \choose M} \zeta_N^{n M}$.  It follows
that 
\[
\sum_{n=1}^N {N \choose M} \zeta_N^{n M} =  {N \choose  M}
\sum_{n=1}^N (\zeta_N^M)^n = 0.
\]
And since this is the coefficient of $x^{N-M}$ in the polynomial 
$\mathbf{v}_1(x) + \mathbf{v}_2(x) + \cdots + \mathbf{v}_N(x)$, we
conclude that
\begin{eqnarray*}
\frac{1}{N} ( \mathbf{v}_1 + \mathbf{v}_2 + \cdots + \mathbf{v}_N )(x)
&=& x^N + 1 \\
&=& \boldsymbol{0}(x),
\end{eqnarray*}
which establishes the lemma.
\end{proof}

\begin{lemma}\label{vertexdistances} 
Let $1 < m, m' < N$.  Then,
$$
\| \mathbf{v}_m - \mathbf{v}_{m'} \|^2 = \| \mathbf{v}_N -
\mathbf{v}_k \|^2 \quad \mbox{where} \quad k = |m' - m|.
$$
Moreover,
$$
\| \mathbf{v}_N - \mathbf{v}_2 \|^2 < \| \mathbf{v}_N - \mathbf{v}_4
\|^2 < \cdots < \| \mathbf{v}_N - \mathbf{v}_{M} \|^2 < 2 {2N
\choose N},
$$
and
$$
\| \mathbf{v}_N - \mathbf{v}_1 \|^2 > \| \mathbf{v}_N - \mathbf{v}_3
\|^2 > \cdots > \| \mathbf{v}_N - \mathbf{v}_{M'} \|^2 >
2 {2N \choose N}, $$ where $M$ is the greatest even integer
not exceeding $N/2$ and $M'$ is the greatest odd integer not exceeding
$N/2$.
\end{lemma}
\begin{proof} 
First we notice that $\| \mathbf{v}_m - \mathbf{v}_{m'} \|^2 =
\|\mathbf{v}_N - \mathbf{v}_k\|^2$ since the isometry 
$R^{N-m}$ takes $\mathbf{v}_m$ to $\mathbf{v}_N$ and $\mathbf{v}_m'$
to $\mathbf{v}_k$.  By (\ref{eq:12}),
\begin{equation}
\label{Parseval result} 
\| \mathbf{v}_N - \mathbf{v}_{k} \|^2 = 
\frac{1}{2 \pi} \int_{-\pi}^{\pi} \left|
\left(e^{i \theta} + 1 \right)^N - \left(e^{i \theta} + \zeta_N^{k}\right)^N \right|^2 \, d\theta 
\end{equation}

Now, if $-\pi < \theta < \pi$ then $\arg(e^{i \theta} + 1) = \theta/2$,
and thus
$$
\arg(e^{i \theta} + \zeta_N^k) = \arg\left(\zeta_N^k \left(1 + e^{i
\theta} \zeta_N^{-k}\right) \right) = \frac{\theta}{2} + \frac{\pi
k}{N}.
$$
It follows that $\arg\left((e^{i \theta} + \zeta_N^k)^N \right) =
(-1)^k \arg\left((e^{i \theta} + 1)^N \right).$ That is, for fixed
$\theta$, $(e^{i \theta} + \zeta_N^k)^N$ and $(e^{i \theta} +
  1)^N$ lie on a line passing through the origin.  From which
we conclude,
\begin{equation}
\label{arg result} \left| \left(e^{i \theta} + 1 \right)^N -
\left(e^{i \theta} + \zeta_N^{k}\right)^N \right|^2 = \left( \left|
\left(e^{i \theta} + 1 \right)^N \right| + (-1)^{k+1} \left|
\left(e^{i \theta} + \zeta_N^{k}\right)^N \right| \right)^2.
\end{equation} Squaring out the right hand side of (\ref{arg result})
and substituting into (\ref{Parseval result}) we find,
\begin{eqnarray} \; \| \mathbf{v}_N - \mathbf{v}_k \|^2 &=& \frac{1}{2
\pi} \int_{-\pi}^{\pi} \left|(e^{i \theta} + 1)^N \right|^2 \, d\theta
+ \frac{1}{2 \pi} \int_{-\pi}^{\pi} \left|(e^{i \theta} + \zeta_N^k)^N
\right|^2 \, d\theta
  \label{split term}\\ && \hspace{1.35cm} + \,  (-1)^{k+1}
\left(\frac{1}{\pi}\int_{-\pi}^{\pi} \left|(e^{i \theta} +
1)^N\right| \left| (e^{i \theta} + \zeta_N^k)^N \right| \, d\theta
\right) \nonumber \\ &=& \label{dist int} 2 {2 N \choose N}
+  (-1)^{k+1} \left(\frac{1}{\pi}\int_{-\pi}^{\pi}
\left|(e^{i \theta} + 1)^N\right| \left| (e^{i \theta} + \zeta_N^k)^N
\right| \, d\theta \right),
\end{eqnarray} where the binomial coefficient in (\ref{dist int})
arises from an application of Parseval's formula to the first two
integrals in (\ref{split term}) together with the familiar formula for
the sum of the squares of the binomial coefficients.

Next we define $f(\theta) = \left|(e^{i \theta} + 1)^N \right|$ and
notice that $f(-\theta) = f(\theta)$ and $f(2\pi k/N - \theta) =
\left|(e^{i \theta} + \zeta_N^k)\right|^N$.  Consequently, we can
write the integral in (\ref{dist int}) as $f\! \ast \!f(2 \pi k/N)$
where
$$
f\! \ast \! f(t) = \frac{1}{\pi} \int_{-\pi}^{\pi} f(\theta) f(t -
\theta) \, d\theta.
$$
The continuity of $f$ implies that $f\! \ast \!f$ is itself
continuous, and likewise the fact that $f$ is even gives $f\! \ast
\!f(-t) = f\!\ast\! f(t)$.  The lemma will be proved by showing that
$f\!\ast\!  f(t)$ is increasing on $(-\pi, 0)$ and decreasing on $(0,
\pi)$.  In order to do this, we use the elementary fact that
$f(\theta) = (2 + 2 \cos{\theta})^{N/2}$ and hence $f$ is increasing
on $(-\pi, 0)$ and decreasing on $(0, \pi)$.  Moreover, $f$ is
differentiable on $(-\pi, \pi)$.  Consequently, we may write
$$
(f \! \ast \! f)'(t) = \frac{1}{\pi} \int_{-\pi}^{\pi} f'(\theta) f(t
- \theta) \, d\theta = \frac{1}{\pi} \int_0^{\pi} f'(\theta) \{ f(t -
\theta) - f(t + \theta) \} \, d\theta,
$$
where the latter equation holds since $f$ is even.  

Notice that for fixed $t$, $f(t - \theta) - f(t + \theta) = 0$ if and
only if $\theta = 0$ or $\theta = \pi$, and $f'(\theta) <
0$ on $(0, \pi)$.  It follows that, for fixed $t$, the integrand never
changes sign.  When $0 < t < \pi$ and $\theta = t/2$ the integrand is
negative.  From which we conclude that $(f \ast f)'(t) < 0$ on $(0,
\pi)$.  Similarly $(f \ast f)'(t) > 0$ on $(-\pi, 0)$, which
establishes the lemma.
\end{proof}

\subsection{The Proof of  Theorem~\ref{thm:3}}

We begin by proving that every isometry of $W_N$ fixes the origin. 

Let $D$ be the distance between $\boldsymbol{0}$ and a vertex of
$W_N$. Let $B$ denote the ball of radius $D$ centred at
$\boldsymbol{0}$ in $\R^{N-1}$, and let $S=\partial B$.  Therefore $S$
is the smallest sphere in $\mathbb{R}^{N-1}$ centred at
$\boldsymbol{0}$ that circumscribes $W_N$.  Notice that by Proposition
\ref{distance prop}, $||\mathbf{w}-\boldsymbol{0}|| \leq D$ for all
$\mathbf{w}\in W_N$ with equality precisely when $\mathbf{w}$ is a
vertex.  Let $T$ be an isometry of $W_N$.  If $T
(\boldsymbol{0})=\mathbf{w}$ for some non-zero $\mathbf{w}$ in $W_N$
then $|| \mathbf{w}-\mathbf{w'}|| \leq D$ for all $\mathbf{w'} \in
W_N$, and in particular, $|| \mathbf{w}- \mathbf{v}_n|| \leq D$ for
$1\leq n\leq N$.  $T(B)$ is the ball of radius $D$ centred at
$\mathbf{w}$, and $W_N \subset B \cap T(B),$ and moreover the vertices
of $W_N$ must lie in $S \cap T(B)$.  And, since $S \cap T(B)$ lies in
the half space $\{ \mathbf{u} \in \R^{N-1} : \mathbf{w} \cdot
\mathbf{u} > 0\}$, we must have $ \mathbf{w} \cdot \mathbf{v}_n > 0$
for $1 \leq n \leq N$.  But this is a contradiction since, by
Lemma~\ref{lemma:3},
\[
0 = \mathbf{w} \cdot \left(\sum_{n=1}^N \mathbf{v}_n\right) =
\sum_{n=1}^N \mathbf{w} \cdot \mathbf{v_n} > 0.
\]
We conclude there is no such $\mathbf{w}$, and that every isometry of
$W_N$ fixes $\boldsymbol{0}$.

Since every isometry fixes the origin, Proposition \ref{distance prop}
implies that any isometry must permute the vertices of $W_N$.  Since
$R$ is transitive on the vertices, composing $T$ with the appropriate
power of $R$ fixes $\mathbf{v}_N$.  By Proposition
\ref{vertexdistances}, $\mathbf{v}_1$ and $\mathbf{v}_{n-1}$ are the
unique vertices farthest from $\mathbf{v}_N$.  Therefore, either
$\mathbf{v}_1$ is fixed or it is exchanged with $\mathbf{v}_{N-1}.$
Recall that $\mathbf{v}_m$ corresponds to $(x+\zeta_N^m)^N$, so
conjugation exchanges $\mathbf{v}_1$ and $\mathbf{v}_{N-1}$.
Composition with $C$, if necessary, results in an isometry, $S$, that
fixes $\mathbf{v}_N$ and $\mathbf{v}_1$.

We will now show that any isometry that fixes both $\mathbf{v}_N$ and
$\mathbf{v}_1$ is the identity.  As $R$ and $C$ are in $D_N$, this
will imply that $T\in D_N$, establishing the result.   As the farthest vertices 
from $\mathbf{v}_1$ are
$\mathbf{v}_N$ and $\mathbf{v}_2$, and $S$ fixes $\mathbf{v}_1$, we conclude that
$S(\mathbf{v}_2)=\mathbf{v}_2$ as $S$ fixes $\mathbf{v}_N$.
Continuing in this fashion, we see that $S$ fixes all vertices.
$S$ fixes the origin, and therefore if we let $\ell_n$ be the 
line segment connecting $\boldsymbol{0}$ and $\mathbf{v}_n$, $S$ fixes $\ell_n$ pointwise.
As a result, $S$ fixes the span of the vertices, which  is $\mathbb{R}^{N-1}$ by 
Lemma \ref{lemma:1}.  We conclude that $S$ is the identity.

\section{The Volume of $W_N$}\label{volume}

\subsection{$\omega$-Conjugate Reciprocal Polynomials}
We will appeal to results and methods from random matrix theory in
order to determine the volume of $W_N$.  If $\lambda_{N-1}$ is
Lebesgue measure on $\R^{N-1}$ then
\begin{equation}
\label{eq:14}
\vol(W_N) = \lambda_{N-1}(W_N) = \int_{W_N} d\lambda_{N-1}(\mathbf{w}).
\end{equation}
Our basic strategy will be to view elements of $W_N$ as polynomials
and employ a change of variables so that we are integrating over the
roots of CR polynomials as opposed to the coefficients.

Along these lines, we must enlarge the set of polynomials under
investigation.  Given a fixed $\omega \in \T$ we say the degree $N$
polynomial $f(x)$ is $\omega$-{\it conjugate reciprocal} (or $\omega$-CR) if
\[
f(x) = \omega x^N \overline{ f(1/\overline{x})}.
\]
In analogy with CR-polynomials, we define the $N -1 \times N-1$ matrix
$X_{N,\omega}$ by
\begin{equation}
\label{eq:7}
X_{N,\omega}[j,k] = \left\{
  \begin{array}{ll} \frac{\sqrt{2}}{2}\left(\delta_{j,k} + \omega
      \delta_{N-j, k}\right) & \mbox{if} \quad 1 \leq j < \lfloor N/2 \rfloor \\ & \\
      \omega^{1/2} \delta_{j,k} & \mbox{if} \quad j = N/2 
      \\ & \\ \frac{\sqrt{2}}{2}\left(i \delta_{N-j,k} - i \omega
        \delta_{j, k}\right) &  \mbox{if} \quad     \lfloor N/2
      \rfloor < j < N. \\ 
\end{array} \right.
\end{equation}
It follows that if $\mathbf{a} \in \R^{N-1}$ and $\mathbf{c} =
X_{N,\omega} \mathbf{a}$, then $(x^N + \omega) + \sum_{n=1}^{N-1} c_n x^{N-n}$
is $\omega$-CR.  We then define $W_{N,\omega}$ to be 
\[
W_{N,\omega} = \bigg\{ \mathbf{a} \in \R^{N-1} : (x^N + \omega) +
\sum_{n=1}^{N-1} c_n x^{N-n} \mbox{ has all roots on \T, where }
\mathbf{c} = X_{N,\omega}\mathbf{a} \bigg\}.
\]
Now consider the map $E_{N,\omega} : \T^{N-1} \rightarrow W_N$ specified by
$\mathbf{a} = E_{N, \omega}(\bs{\xi}) = X_{N,\omega}^{-1} \mathbf{c}$ where
$\mathbf{c}$ is obtained from $\boldsymbol{\xi}$ by
\begin{equation}
\label{eq:9}
\left(x - \frac{(-1)^N \omega}{\xi_1 \xi_2 \cdots \xi_{N-1}}\right)
\prod_{n=1}^{N-1} (x - \xi_n) = (x^N + \omega) + \sum_{n=1}^{N-1} c_n
x^{N-n}.
\end{equation}
In order to determine the volume of $W_N$ we need to compute the
(absolute value of the) Jacobian of the map $E_{N,1}.$  It shall be
convenient, and no more difficult, to compute the Jacobian of
$E_{N,\omega}$ for arbitrary $\omega$.
\begin{lemma}
Let $\omega \in \T$.  The absolute value of the Jacobian of $E_{N, \omega}$
is given by
\begin{equation}
\label{eq:10}
| \Jac E_{N,\omega} (\xi_1, \xi_2, \ldots, \xi_{N-1}) | = 
 \prod_{1 \leq m < n \leq N} | \xi_n - \xi_m |,
\end{equation}
where
\begin{equation}
\label{eq:5}
\xi_N = \frac{(-1)^N \omega}{\xi_1 \xi_2 \cdots \xi_{N-1}}.
\end{equation}
Moreover,
\begin{equation}
\label{eq:11}
| \Jac E_{N,\omega} (\xi_1, \xi_2, \ldots, \xi_{N-1}) | = | \Jac E_{N,1} (\xi_1, \xi_2, \ldots, \xi_{N-1}) |.
\end{equation}
\end{lemma}
\begin{proof}
By (\ref{eq:9}),
\begin{eqnarray}
c_n &=& (-1)^{n} e_n(\xi_1, \xi_2, \ldots, \xi_N) \nonumber \\
&=& (-1)^{n} \big(e_n(\xi_1, \xi_2, \ldots, \xi_{N-1}) + \xi_N
e_{n-1}(\xi_1, \xi_2, \ldots, \xi_{N-1}) \big), \label{eq:4}
\end{eqnarray}
where $e_n$ is the $n$th elementary symmetric function.  The Jacobian
of $E_{N, \omega}$ is given by 
\[
\Jac E_{N, \omega} = \det \left[ \frac{\partial a_n}{\partial \xi_m}
\right]_{n,m=1}^{N-1} = \det X_{N,\omega}^{-1} \det \left[ \frac{\partial
 c_n}{\partial \xi_m} \right]_{n,m=1}^{N-1}.
\]
From (\ref{eq:4}) an easy calculation reveals
\[
\frac{\partial c_n}{\partial \xi_m} = (-1)^{n} \left(1 -
  \frac{\xi_N}{\xi_m}\right) e_{n-1,m}
\]
where $e_{n-1,m} = e_{n-1,m}(\xi_1, \xi_2, \ldots, \xi_{N-1})$ is the
$(n-1)$st elementary symmetric function in all variables except $\xi_m$,
and we use the convention that $e_{0,m} = 1$.  We conclude that
\begin{equation}
\label{eq:13}
\Jac E_{N, \omega} = \det X_{N,\omega}^{-1} \left(\prod_{\ell=1}^{N-1}
\xi_{\ell}^{-1} (\xi_N - \xi_{\ell})\right) \det \left[(-1)^{n-1}
e_{n-1,m} \right]_{n,m=1}^{N-1}.
\end{equation}
Let $U$ be the $N-1 \times N-1$ matrix
$\left[(-1)^{n-1} e_{n-1,m}\right]_{n,m=1}^{N-1}$ and notice that
the $m$th column of $U$ is comprised of the coefficients of the
polynomial
\[
f_m(x) = \prod_{\ell=1 \atop \ell \neq m}^{N-1} (x - \xi_\ell) =
\sum_{n=1}^{N-1} (-1)^{n-1} e_{n-1,m} x^{N-1-n}
\]
and clearly
\[
f_m(\xi_{k}) = \delta_{m,k} \prod_{\ell=1 \atop \ell \neq
 m}^{N-1} (\xi_{k} - \xi_{\ell}).
\]
Thus, if we set $V = [ \xi_m^j ]_{m,j=1}^{N-1}$ (that is $V$ is the
Vandermonde matrix in the variables $\xi_1, \xi_2, \ldots,
\xi_{N-1}$), then $VU$ is a diagonal matrix, and
\[
\det VU = (-1)^{N-1 \choose 2} \prod_{1 \leq m < n < N} (\xi_n - \xi_m)^2.
\]
Using the familiar formula for the Vandermonde determinant we see
\[
\det U =  (-1)^{N-1 \choose 2} \prod_{1 \leq m < n < N} (\xi_n - \xi_m).
\]
It is easy to verify from (\ref{eq:8}) and (\ref{eq:7}) that $\det
X_{N,\omega} = \omega^{N/2} \det X_N$ and hence (\ref{eq:13}) leads us to
\[
|\Jac E_{N,\omega}| = \prod_{1 \leq m < n \leq N} |\xi_n - \xi_m|.
\]

To prove (\ref{eq:11}) we note that
\[
(-1)^N \prod_{n=1}^N \omega^{-1/N} \xi_n = \omega^{-1} \left( (-1)^N
  \prod_{n=1}^N \xi_n \right) = \omega^{-1} \omega = 1.
\]
Thus we may transform an $\omega$-CR polynomial into a CR polynomial by
multiplying each of the roots by $\omega^{-1/N}$.  But this corresponds
to multiplying $\Jac E_{N,\omega}$ by a complex number of modulus 1.
\end{proof}

\subsection{The Proof of Theorem ~\ref{thm:4}}
\label{sec:proof-theor-refthm:4}

\sloppypar{
Given $\mathbf{w}(x) = \prod_{n=1}^N (x - \xi_n) \in \overcirc{W_N}$
there are $N!$ different choices of root vectors $\mathbf{\xi} \in
\T^{N-1}$ associated to $\mathbf{w}$.  That 
is, since one of the roots is determined by the others, there are ${N
  \choose N-1}$ ways of choosing $N-1$ independent roots, and then
$(N-1)!$ ways of ordering them.  Setting $\xi_n = e^{i \theta_n}$, it
follows from (\ref{eq:14}) that
\begin{eqnarray*}
\vol(W_N) &=& \frac{1}{N!} \int_{0}^{2 \pi} \cdots \int_{0}^{2 \pi}
\left| \Jac E_{N,1}(e^{i \theta_1}, e^{i \theta_2},
  \ldots, e^{i \theta_{N-1}} ) \right| d\theta_1 \cdots d\theta_{N-1} \\
&=& \frac{1}{N!} \int_0^{2\pi} \cdots \int_{0}^{2 \pi} 
\left( \frac{1}{2\pi} \int_{0}^{2 \pi} \left| \Jac E_{N,1}(e^{i
      \theta_1}, e^{i \theta_2},   \ldots, e^{i \theta_{N-1}} ) \right|
  d\vartheta \right) d\theta_1 \cdots  d\theta_{N-1}.
\end{eqnarray*}
By (\ref{eq:11}) we may replace $\left| \Jac
E_{N,1} \right|$ with $\left| \Jac E_{N,e^{i\vartheta}} \right|$
\[
\vol(W_N) = \frac{1}{N!} \int_0^{2\pi} \cdots \int_{0}^{2 \pi} 
\left( \frac{1}{2\pi} \int_{0}^{2 \pi} \left| \Jac E_{N,e^{i \vartheta}}(e^{i
      \theta_1}, e^{i \theta_2},   \ldots, e^{i \theta_{N-1}} ) \right|
  d\vartheta \right) d\theta_1 \cdots  d\theta_{N-1}.
\]
By the change of variables $\theta_N = N\pi + \vartheta - (\theta_1 + \cdots + \theta_{N-1})$ and (\ref{eq:10}) we discover that
\[
\vol(W_N) = \frac{1}{2 \pi N!} \int_{0}^{2\pi}
\int_{0}^{2\pi} \cdots \int_{0}^{2\pi} \prod_{1 \leq m < n \leq N} |
e^{i \theta_n} - e^{i \theta_m} | \, d\theta_1 \, d\theta_2 \cdots
d\theta_N.
\]
The value of this integral has been calculated by F.~Dyson
in the context of random matrix theory \cite{MR0278668}.  Using
Dyson's value for this integral, we have
\[
\vol(W_N) = \frac{2^{N-1} \pi^{(N-1)/2}}{\Gamma(\frac{N+1}{2})},
\]
which is the volume of the $N-1$ dimensional ball of radius $2$.
}

\section{Acknowledgments} 

The authors would like to thank David Farmer for suggesting the
problem of finding the volume of $W_N$, and for numerous valuable
discussions regarding CR polynomials.

\bibliography{bibliography}

\vspace{.25cm}
\noindent\rule{4cm}{.5pt}

\vspace{.25cm} {\small \noindent 
{\sc Kathleen L. Petersen} \\
{{ Queen's University} \\
 { Kingston, Ontario} \\
 email: {\tt petersen@mast.queensu.ca}}

\vspace{.5cm} {\small \noindent 
{\sc Christopher D. Sinclair} \\
{{ Pacific Institute for the Mathematical Sciences} \\
 { Vancouver, British Columbia} \\
 email: {\tt sinclair@math.ubc.ca}}

\end{document}

%% file: preamble.tex
\usepackage{amsmath,amsthm,amsfonts,amscd,mathrsfs,pifont} 
\usepackage{eucal}  
\usepackage{verbatim}     
\usepackage[dvips]{graphicx}
\usepackage{hyperref}

\newtheorem{thm}{Theorem}[section]
\newtheorem{cor}[thm]{Corollary}

\newtheorem{lemma}[thm]{Lemma}
\newtheorem{prop}[thm]{Proposition}

\theoremstyle{definition}

\theoremstyle{remark}

\newcommand{\BB}[1]{\ensuremath{\mathbb{#1}}}

\newcommand{\R}{\ensuremath{\BB{R}}}
\newcommand{\Z}{\ensuremath{\BB{Z}}}
\newcommand{\C}{\ensuremath{\BB{C}}}
\newcommand{\T}{\ensuremath{\BB{T}}}

\newcommand{\bs}{\ensuremath{\boldsymbol}}

\DeclareMathOperator{\vol}{vol}

\DeclareMathOperator{\Isom}{Isom}
\DeclareMathOperator{\Jac}{Jac}
\DeclareMathOperator{\disc}{disc}